\documentclass[english]{article}
\usepackage{amsfonts}
\usepackage[latin9]{inputenc}
\usepackage{amsmath}
\usepackage{amssymb}
\usepackage{url}
\usepackage{graphicx}
\usepackage[pdfstartview=FitH]{hyperref}
\usepackage{amsthm}
\usepackage{authblk}


\usepackage{babel}

\begin{document}
\title{Vanishing of cohomology and property (T) for groups acting on weighted simplicial complexes}
\author{Izhar Oppenheim}
\affil{Department of Mathematics\\
 The Ohio State University  \\
 Columbus, OH 43210, USA \\
E-mail: izharo@gmail.com}

\newtheorem{theorem}{Theorem}[section]
\newtheorem{lemma}[theorem]{Lemma}
\newtheorem{definition}[theorem]{Definition}
\newtheorem{claim}[theorem]{Claim}
\newtheorem{example}[theorem]{Example}
\newtheorem{remark}[theorem]{Remark}
\newtheorem{proposition}[theorem]{Proposition}
\newtheorem{corollary}[theorem]{Corollary}

\maketitle

\begin{abstract}
We extend Ballmann and \'{S}wiatkowski's work on cohomology of groups acting on simplicial complexes and provide further vanishing results of cohomologies and of property (T). In particular, we give a new criterion for property (T) for groups acting on an $n$-dimensional simplcial complex.  
\end{abstract}

\section{Introduction}

For a finite graph $L$ with a set of vertices $V_{L}$, the Laplacian
of the graph $\Delta^+$ is an operator on the space of real valued
functions on $V_{L}$ which is defined as \[
\Delta^+ f(v)=f(v)-\frac{1}{m(v)}\sum_{u\sim v}f(u)\]
where $m(v)$ is the valance of $v$ and $u\sim v$ means that there
is an edge connecting $u$ and $v$. The Laplacian is a positive operator and we denote by $\lambda(L)$ its smallest positive eigenvalue. One can generalize the definition of the Laplacian in two ways - first, one can put weights on the edges and second, one can generalize the definition so the Laplacian will be defined for a simplicial complex  $X$ of any dimension. For such a complex the (weighted) Laplacian is again a positive operator and we denote by $\lambda(X)$ its smallest positive eigenvalue.

Ballmann and \'{S}wiatkowski in \cite{BS} and independently \.{Z}uk's in \cite{Zuk} gave criteria
for the vanishing of cohomologies of a group $\Gamma$ acting on a simplicial complex $X$, by considering the values of $\lambda$ for the links of $X$ (in both cases, the authors were relaying on the previous work of Garland \cite{Gar}). \\
In \cite{DJ1} the above results were generalized by Dymara and Januszkiewicz to a more general setting in which $\Gamma$ isn't necessarily discrete but just locally compact and unimodular. \\
The case which is maybe the most interesting is the vanishing of the first cohomology, because this is equivalent to property (T). \\
In this paper, we shall generalize the criterion given in \cite{BS} to a criterion for weighted Laplacian and in this slightly more general setting we will prove new vanishing of cohomology results (including the fulfilment of property (T)). Our results are new not only in the weighted setting, but are also in the setting given in \cite{BS}. Namely (in the setting of \cite{BS}), we will prove the following theorem:

\begin{theorem}
\label{introthm}
Let $X$ be a contractible pure simplicial complex of dimension $n$, such that all the links of $X$ of dimension $>0$ are connected, and such that there is an integer $0 \leq r \leq n-2$ such that the link of every $r$-simplex is finite. 
Let $\Gamma$ be a locally compact, unimodular group acting cocompactly by
automorphisms on $X$. Assume that for every simplex $\eta$ of dimension $r$ the group $\Gamma_\eta$ is compact ($\Gamma_\eta \subset \Gamma$ is the point-wise stabilizer of $\eta$ in $\Gamma$). For a simplex $\tau$ of dimension $k-1 \geq r$, denote by $X_{\tau}$ the link of $\tau$ which is a simplicial complex and denote by $\lambda(X_{\tau})$ the smallest positive eigenvalue of the (untwisted) Laplacian on $X_{\tau}$. \\
Let $r < k \leq n-1$, if $\lambda(X_{\tau}) > \frac{k(n-k)}{k+1}$ for every simplex $\tau$ of dimension $k-1 \geq r$,
then:
\begin{enumerate}
\item For every $r < j \leq k$ we have $H^{j}(X,\rho)=0$ for any unitary representation $\rho$ of $\Gamma$.
\item $\Gamma$ has property (T).
\end{enumerate}

\end{theorem}

The power of this theorem is that it allows us to prove property (T) just by calculating the Laplacian's smallest positive eigenvalue for links of simplices of dimension $n-2$ (in those cases the links are graphs, so the eigenvalues are easier to compute). Also, in the case that $X$ is locally finite we have that $H^{j}(X,\rho)=H^{j}(\Gamma,\rho)$ (see remark below) and therefore in that case we have the following corollary:
\begin{corollary}
Under the assumptions of theorem \ref{introthm} if $X$ is locally finite and the action of $\Gamma$ is proper (i.e. $r=0$ in the theorem) we have the following: \\
Let $0 < k \leq n-1$, if $\lambda(X_{\tau}) > \frac{k(n-k)}{k+1}$ for every simplex $\tau$ of dimension $k-1 \geq 0$,
then for every $0 < j \leq k$ we have $H^{j}(\Gamma,\rho)=H^{j}(X,\rho)=0$ for any unitary representation $\rho$ of $\Gamma$.
\end{corollary}

\begin{remark}
In general, the connection between $H^* (\Gamma, \rho)$ and $H^* (X, \rho)$ is given by a spectral sequence converging $H^* (\Gamma, \rho)$ with $E_1^{p,q} =\oplus_{\tau \in \Sigma (p, \Gamma)} H^q (\Gamma_\tau, \rho) $ (where $\Sigma (p, \Gamma)$ are $p$-simplices in $X/\Gamma$) and $E_2^{p,0} = H^p (X,\rho)$ (see \cite{DJ1}[Lemma 2.1] and also \cite{DJ2}[Section 3] for a more complete discussion). In the case where $X$ is locally finite we get that $G_\tau$ are always compact and therefore  $H^q (\Gamma_\tau, \rho) =0$ and $H^j (X, \rho) = H^j (\Gamma, \rho)$ for every $1 \leq j \leq n-1$. In the non locally finite case, we can't assume $H^q (\Gamma_\tau, \rho) =0$ and therefore it is not enough to show $H^j (X, \rho) =0$ to show vanishing of $H^j (\Gamma, \rho)$. This is resolved in this article for the case $j=1$ (property (T)), by replacing $X$, which isn't locally finite, by $X'$, which is locally finite and has a similar combinatorial structure (i.e. the spectral gap of the Laplacian eigenvalues is large enough to ensure vanishing of the first cohomology). The "price" one pays when replacing $X$ by $X'$ is adding higher cohomologies - even if $X$ is assumed contractible, $X'$ is in general only connected and simply connected (in general $X'$ is likely to have non trivial cohomology groups $H^k (X')$ for $k>1$).
\end{remark}
\begin{remark}
Theorems similar to the one stated above were proven in several articles (other then \cite{BS} and \cite{Zuk} that were already mentioned).  We will not give a detailed account about all those articles, but just mention a few for comparison:
\begin{enumerate}
\item Garland \cite{Gar} considered the case of a group $\Gamma$ acting on a building whose all of its vertices has finite links. In this case, \cite{Gar} showed that for every $0 \leq k \leq n-1$ there is a $\varepsilon >0$ such that if all the positive Laplacian eigenvalues of links of dimension $1$ are strictly larger than $1- \varepsilon$ then $H^{k}(\Gamma,\rho)=0$. But the nature of this result was less quantitative and more qualitative: there is no sharp estimation of $\varepsilon$.
\item Dymara and Januszkiewicz \cite{DJ2} considered the case of a group acting on a simplicial complex where the fundamental domain is a single simplex and there is an integer $0 \leq r \leq n-2$ such that the link of every $r$-simplex is finite. In this case they proved that if the smallest positive Laplacian eigenvalue of links of dimension $1$ are large enough, then the following  vanishing results hold: $H^{j}(X,\rho)=0$ for every $0< j \leq n-1$ (\cite{DJ2}[theorem 5.1]) and $H^{j}(\Gamma,\rho)=0$ for $0< j \leq n-r-1$  (\cite{DJ2}[theorem 5.2]). The reader should note that although this result is more powerful that the results of this article in many aspects, it requires two stronger assumptions: first, the fundamental domain has to be a single simplex. Second, the spectral gap needed is quite large: the smallest positive Laplacian eigenvalue should be strictly larger than $1- \frac{13}{28^n}$ to ensure the vanishing of the cohomologies.       
\item Kassabov \cite{Kas} considered a setting similar to the one considered in \cite{DJ2}. In that setting, \cite{Kas} proves that if the smallest positive eigenvalue is larger than $\frac{n-1}{n}$, then the group has property (T). In fact, \cite{Kas} proves a little more general criterion which takes into account the interplay between different $1$-dimensional links. However, the proof of \cite{Kas} relays on the fact that the fundamental domain is a single simplex and we do not think that the method used by Kassabov can be generalized to the setting of this article, where the assumption is that the fundamental domain is just compact (but can contain more than a single simplex).  
\end{enumerate}
\end{remark}

\textbf{Structure of the paper.} Section 2 is devoted to introducing the framework developed in \cite{BS} and \cite{DJ1} with some generalizations (weights) and additions (restriction). In section 3, we prove vanishing results of $H^{j}(X,\rho)$ given a lower bound on the Laplacian positive eigenvalues (as in the first implication of theorem \ref{introthm}). In section 4, we prove that the criterion given in section 3 for the vanishing of $H^{j}(X,\rho)$ also implies property (T) even in the case where $X$ is not locally finite (as in the section implication of theorem \ref{introthm}).  \\ 

\textbf{Acknowledgement.} The results presented in this article were part of the author Ph.D. thesis written under the supervision of prof. Uri Bader at the Technion - Israel Institute of Technology. The author would also like to thank the anonymous referee for his many useful suggestions and for pointing out errors in the original manuscript.

\section{Framework}

Here we introduce a slight generalization of the framework constructed in \cite{BS} and \cite{DJ1}. We should note that this generalization was already considered by Wang in \cite{Wang} but for completeness we will prove all the propositions that differ (by a constant) from those of \cite{BS} and \cite{DJ1}. \\ 
Throughout this paper $X$ is a $n$-dimensional simplicial complex with the following properties:
\begin{enumerate}
\item $X$ is pure $n$-dimensional, i.e., every simplex in $X$ is contained in at least one $n$-dimensional simplex.
\item X is contractible.
\item all the links of $X$ of dimension $>0$ are connected.
\item There is a constant integer $ 0 \leq r \leq n-2$ such that all the links of simplices of dimension $r$ are finite.
\item $X$ is has a weight function $m$ - see definition \ref{weightedDef} below.
\end{enumerate}
 Also, throughout this paper, $\Gamma$ is a locally compact, unimodular group of automorphisms of $X$ acting cocompactly on $X$ such that for every simplex $\sigma$ of dimension $r$ the stabilizer group of $\sigma$, denoted by $\Gamma_\sigma$, is compact. Last, $\rho$ is a unitary representation of $\Gamma$ on a complex Hilbert space $H$.

\subsection{general settings}

Following \cite{BS} we introduce the following notations:
\begin{enumerate}
\item For $0\leq k\leq n$, denote by $\Sigma(k)$ the set of ordered $k$-simplices (i.e. $\sigma \in \Sigma(k)$ is an ordered $(k+1)$-tuple of vertices that form a $k$-simplex in $X$) and choose a set $\Sigma(k,\Gamma)\subseteq\Sigma(k)$ of representatives
of $\Gamma$-orbits.
\item For a simplex $\sigma\in\Sigma(k)$, denote by $\Gamma_{\sigma}$
the stabilizer of $\sigma$ and by $\vert \Gamma_{\sigma} \vert$ the Haar measure of $\Gamma_{\sigma}$.  
\item For $0\leq k\leq n$, denote by $C^{k}(X,\rho)$ the space of simplicial
$k$-cochains of $X$ which are twisted by $\rho$, that is, $\phi\in C^{k}(X,\rho)$
is an alternating map on ordered $k$-simplices of $X$ with values
in $H$ such that \[
\forall\gamma\in\Gamma,\forall x\in\Sigma(k),\rho(\gamma)\phi(x)=\phi(\gamma x)\]

\end{enumerate}

The following proposition is taken from \cite{BS}, \cite{DJ1}:

\begin{proposition}
\cite[Lemma 1.3]{BS}, \cite[Lemma 3.3]{DJ1} For $r \leq l<k\leq n$, let $f=f(\tau,\sigma)$ be
a $\Gamma$-invariant function on the set of pairs $\left(\tau,\sigma\right)$,
where $\tau$ is an ordered $l$-simplex and $\sigma$ is an ordered
$k$-simplex with $\tau\subset\sigma$ Then \[
\sum_{\sigma\in\Sigma(k,\Gamma)}\sum_{\begin{array}{c}
{\scriptstyle \tau\in\Sigma(l)}\\
{\scriptstyle \tau\subset\sigma}\end{array}}\frac{f(\tau,\sigma)}{\left|\Gamma_{\sigma}\right|}=\sum_{\tau\in\Sigma(l,\Gamma)}\sum_{\begin{array}{c}
{\scriptstyle \sigma\in\Sigma(k)}\\
{\scriptstyle \tau\subset\sigma}\end{array}}\frac{f(\tau,\sigma)}{\left|\Gamma_{\tau}\right|}\]

\end{proposition}

The reader should note, that from now on we will use the above proposition to change the order of summation without mentioning it explicitly.

\begin{definition}
\label{weightedDef}
A simplicial complex $X$ is called weighted if there is an equivariant function $m : \bigcup_{k \geq r} \Sigma(k) \rightarrow \mathbb{R}^+$ (called the weight function) such that:
\begin{enumerate}
\item For every $\tau = (v_0,...,v_k)$ and for every permutation $\sigma \in S_{k+1}$ we have $m ((v_0,...,v_k)) = m ((v_{\sigma (0)},...,v_{\sigma (k)})$.
\item For every $r \leq k \leq n-1$ there is a constant $C_k$ such that for every $\tau \in \Sigma (k)$ we have the following equality
$$ \sum_{\sigma \in \Sigma (k+1), \tau \subset \sigma} m( \sigma ) = (k+2)! C_k m (\tau )$$
Where $\tau \subset \sigma$ means that all the vertices of $\tau$ are contained in $\sigma$ (with no regard to the ordering). 
\end{enumerate} 
\end{definition}

\begin{example}
\label{BSweights}
In \cite{BS} the function $m$ was defined as: for every $\tau \in \Sigma (k)$, $m (\tau)$ is the number of (unordered) simplices of dimension $n$ that contain $\tau$. In that case, $m$ is the constant function $1$ on $\Sigma (n)$ and $C_k = n-k$.   
\end{example}

\begin{remark}
There is a lot of freedom in our definition of the weight function. Without loss of generality, one can always normalize the weight function such that $C_k = 1$ for all $k$ (and this is the setting used in \cite{Wang}). It is obvious that in the normalized case the function $m$ is determined by its values on $\Sigma (n)$. We chose not to normalize the weight function in this paper as a matter of convenience and so that the reader could easily compare our results to those of \cite{BS}.  
\end{remark}

For $k \geq r$ define an Hermitian form on $C^{k}(X,\rho)$ as
\[
\left\langle \phi,\psi\right\rangle :=\sum_{\sigma\in\Sigma(k,\Gamma)}\frac{m(\sigma)}{(k+1)!\left|\Gamma_{\sigma}\right|}\left\langle \phi(\sigma),\psi(\sigma)\right\rangle \]
With this Hermitian form $C^{k}(X,\rho)$  is a complex Hilbert space. \\
To distinguish the norm of $H$ from the norm of $C^{k}(X,\rho)$
we will use $\left|.\right|$ to denote the norm of $H$ (i.e. $\left\langle \phi(\sigma),\phi(\sigma)\right\rangle =\left|\phi(\sigma)\right|^{2}$) and $\Vert . \Vert$ to denote the norm of $C^{k}(X,\rho)$. \\
For $r \leq k <n$, the \emph{differential} $d:C^{k}(X,\rho)\rightarrow C^{k+1}(X,\rho)$
is given by \[
d\phi(\sigma):=\sum_{i=0}^{k+1}(-1)^{i}\phi(\sigma_{i}),\;\sigma\in\Sigma(k+1)\]
Where $\sigma_{i}=\left(v_{0},...,\hat{v_{i}},...,v_{k+1}\right)$
for $\left(v_{0},...,v_{k+1}\right)=\sigma\in\Sigma(k+1)$. \\
Denote by $\delta$ the adjoint operator of $d$: 
$$\delta : C^{k+1}(X,\rho)\rightarrow C^{k}(X,\rho)$$
Also, denote
$$\Delta^+ = \delta d : C^k (X, \rho) \rightarrow C^k (X, \rho)$$
Define the $k$-th cohomology of $X$ twisted by $\rho$ as
$$ H^k (X, \rho ) = ker (d: C^k (X, \rho) \rightarrow C^{k+1} (X, \rho)) / im (d: C^{k-1} (X, \rho) \rightarrow C^{k} (X, \rho))$$

\begin{proposition}
For $X$,$\Gamma$, $\rho$ as above and $k \geq r$ we have:
\begin{enumerate}
\item (equivalent to \cite[Proposition 1.5]{BS}) The differential it is a bounded operator.
\item (equivalent to \cite[Proposition 1.6]{BS}) The adjoint operator of $d$, denoted
by $\delta:C^{k+1}(X,\rho)\rightarrow C^{k}(X,\rho)$ is \[
\delta\phi(\tau)=\sum_{\begin{array}{c}
{\scriptstyle v \in\Sigma(0)}\\
{\scriptstyle v \tau\in \Sigma (k+1)}\end{array}}\frac{m(v\tau)}{m(\tau)}\phi(v\tau),\;\tau\in\Sigma(k)\]
where $v\tau=(v,v_{0},...,v_{k})$ for $\tau=(v_{0},...,v_{k})$ 
\item (equivalent to \cite[Corollary 1.7]{BS}) For $\phi\in C^{k}(X,\rho)$ and $\sigma\in\Sigma(k)$,
\[
\Delta^+ \phi (\sigma) = \delta d\phi(\sigma)=C_k \phi(\sigma)-\sum_{\begin{array}{c}
{\scriptstyle v \in\Sigma(0)}\\
{\scriptstyle v \sigma \in \Sigma (k+1)}\end{array} }\sum_{0\leq i\leq k}(-1)^{i}\frac{m(v\sigma)}{m(\sigma)}\phi(v\sigma_{i})\]

\end{enumerate}
\end{proposition}

\begin{proof}
\begin{enumerate}
\item For every $\phi \in C^k (X, \rho)$ we have
$$ \Vert d \phi \Vert^2 = \sum_{\sigma \in \Sigma (k+1, \Gamma)} \dfrac{m(\sigma )}{(k+2)! \vert \Gamma_\sigma \vert} \vert \sum_{i=0}^{k+1}(-1)^{i}\phi(\sigma_{i}) \vert^2 \leq $$
$$ \leq \sum_{\sigma \in \Sigma (k+1, \Gamma)} \dfrac{m(\sigma )}{(k+2)! \vert \Gamma_\sigma \vert} (k+2) \sum_{i=0}^{k+1} \vert \phi(\sigma_{i}) \vert^2 \leq$$
$$ \leq \sum_{\sigma \in \Sigma (k+1, \Gamma)} \dfrac{(k+2) m(\sigma )}{(k+2)! (k+1)! \vert \Gamma_\sigma \vert}  \sum_{ \tau \in \Sigma (k), \tau \subset \sigma} \vert \phi(\tau) \vert^2 \leq$$ 
$$ \leq \sum_{\tau \in \Sigma (k, \Gamma)} \dfrac{(k+2) \vert \phi(\tau) \vert^2 }{(k+2)! (k+1)! \vert \Gamma_\tau \vert}  \sum_{ \sigma \in \Sigma (k+1), \tau \subset \sigma} m(\sigma ) \leq$$  
$$ \leq (k+2) C_k \sum_{\tau \in \Sigma (k, \Gamma)} \dfrac{m (\tau ) \vert \phi(\tau) \vert^2 }{ (k+1)! \vert \Gamma_\tau \vert}  = (k+2) C_k \Vert \phi \Vert^2$$  

\item For $\sigma \in \Sigma (k+1)$ and $\tau \subset \sigma, \tau \in \Sigma (k)$ denote by $[\sigma : \tau ]$ the incidence coefficient of $\tau$ with respect to $\sigma$, i.e. if $\sigma_i$ has the same vertices as $\tau$ then for every $\psi \in C^k (X, \rho)$ we have $[ \sigma : \tau ] \psi (\tau ) = (-1)^i \psi (\sigma_i)$. Take $\phi \in C^{k+1} (X, \rho)$ and $\psi \in C^k (X, \rho)$, then we have
$$ \langle d \psi , \phi \rangle = \sum_{\sigma \in \Sigma (k+1, \Gamma)} \dfrac{m(\sigma )}{(k+2)! \vert \Gamma_\sigma \vert} \langle \sum_{i=0}^{k+1}(-1)^{i}\psi(\sigma_{i}) , \phi ( \sigma ) \rangle = $$
$$ = \sum_{\sigma \in \Sigma (k+1, \Gamma)} \dfrac{m(\sigma )}{(k+1)! (k+2)! \vert \Gamma_\sigma \vert} \langle \sum_{\begin{array}{c}
{\scriptstyle \tau \in\Sigma(k)}\\
{\scriptstyle \tau \subset \sigma }\end{array} } [\sigma : \tau] \psi(\tau) , \phi ( \sigma ) \rangle = $$
$$ = \sum_{\sigma \in \Sigma (k+1, \Gamma)} \dfrac{m(\tau )}{(k+1)! \vert \Gamma_\sigma \vert} \sum_{\begin{array}{c}
{\scriptstyle \tau \in\Sigma(k)}\\
{\scriptstyle \tau \subset \sigma }\end{array}} \langle  \psi(\tau) ,\dfrac{[\sigma : \tau] m(\sigma )}{m(\tau ) (k+2)!} \phi ( \sigma ) \rangle = $$
$$ = \sum_{\tau \in \Sigma (k, \Gamma)} \dfrac{m(\tau )}{(k+1)! \vert \Gamma_\tau \vert} \sum_{\begin{array}{c}
{\scriptstyle \sigma \in\Sigma(k+1)}\\
{\scriptstyle \tau \subset \sigma }\end{array}} \langle  \psi(\tau) ,\dfrac{[\sigma : \tau] m(\sigma )}{m(\tau ) (k+2)!} \phi ( \sigma ) \rangle = $$
$$ = \sum_{\tau \in \Sigma (k, \Gamma)} \dfrac{m(\tau )}{(k+1)! \vert \Gamma_\tau \vert}  \langle  \psi(\tau) ,\sum_{\begin{array}{c}
{\scriptstyle \sigma \in\Sigma(k+1)}\\
{\scriptstyle \tau \subset \sigma }\end{array}} \dfrac{[\sigma : \tau] m(\sigma )}{m(\tau ) (k+2)!} \phi ( \sigma ) \rangle = $$
$$ = \sum_{\tau \in \Sigma (k, \Gamma)} \dfrac{m(\tau )}{(k+1)! \vert \Gamma_\tau \vert}  \langle  \psi(\tau) ,\sum_{\begin{array}{c}
{\scriptstyle v \in\Sigma(0)}\\
{\scriptstyle v\tau\in\Sigma(k+1) }\end{array} }\frac{m(v\tau)}{m(\tau)}\phi(v\tau) \rangle $$

\item For every $\phi \in C^k (X, \rho)$ and every $\sigma \in \Sigma (k)$ we have:
$$ \delta d \phi (\sigma ) = \sum_{\begin{array}{c}
{\scriptstyle v \in\Sigma(0)}\\
{\scriptstyle v\sigma \in\Sigma(k+1) }\end{array} }\frac{m(v\sigma)}{m(\sigma)} d \phi(v\sigma) = $$ 
$$ = \sum_{\begin{array}{c}
{\scriptstyle v \in\Sigma(0)}\\
{\scriptstyle v\sigma \in\Sigma(k+1) }\end{array}}\frac{m(v\sigma)}{m(\sigma)} \phi (\sigma ) - \sum_{\begin{array}{c}
{\scriptstyle v \in\Sigma(0)}\\
{\scriptstyle v\sigma \in\Sigma(k+1) }\end{array} }\sum_{0\leq i\leq k}(-1)^{i}\frac{m(v\sigma)}{m(\sigma)}\phi(v\sigma_{i}) =$$
$$ = \sum_{\begin{array}{c}
{\scriptstyle \gamma \in\Sigma(k+1)}\\
{\scriptstyle \sigma \subset \gamma }\end{array}}\frac{m(\gamma)}{(k+2)! m(\sigma)} \phi (\sigma ) - \sum_{\begin{array}{c}
{\scriptstyle v \in\Sigma(0)}\\
{\scriptstyle v\sigma \in\Sigma(k+1) }\end{array}}\sum_{0\leq i\leq k}(-1)^{i}\frac{m(v\sigma)}{m(\sigma)}\phi(v\sigma_{i}) =$$
$$ = C_k \phi (\sigma ) - \sum_{\begin{array}{c}
{\scriptstyle v \in\Sigma(0)}\\
{\scriptstyle v\sigma \in\Sigma(k+1) }\end{array}}\sum_{0\leq i\leq k}(-1)^{i}\frac{m(v\sigma)}{m(\sigma)}\phi(v\sigma_{i})$$
\end{enumerate}
\end{proof}

\begin{theorem}
\label{CohomologyVanish}
Let $X$, $\Gamma$ and $\rho$ be as above and let $k>r$. If there is an $\varepsilon >0$ such that for every $\phi \in C^k (X,\rho), d \phi =0$ one has $\Vert \delta \phi \Vert^2 \geq \varepsilon \Vert \phi \Vert^2$ then $H^k (X,\rho)=0$. If in addition $r=0$ (i.e., $X$ is locally finite and the action of $\Gamma$ is proper) then $H^k (\Gamma, \rho) = H^k (X,\rho)=0$.
\end{theorem}

\begin{proof}
This is a general theorem that doesn't require any adaptation to the weighted setting so we omit its proof. Instead, we refer the reader to \cite{BS}[section 2] and \cite{DJ1}[Theorem 3.1]  for proof of the fact that $\Vert \delta \phi \Vert^2 \geq \varepsilon \Vert \phi \Vert^2$ for every $\phi \in C^k (X,\rho), d \phi =0$ implies $H^k (X,\rho)=0$ and to \cite{DJ1}[Lemma 2.1] for proof that in the locally finite case (when the group action is proper) one has $H^k (\Gamma,\rho)=H^k (X,\rho)$.
\end{proof}

\begin{proposition}
\label{ConstantEqui}
Let $X$, $\Gamma$, $\rho$ as above. Assume further that $r=0$, i.e. $X$ is locally finite and $\Gamma$ acts properly on $X$. For $\phi \in C^0 (X,\rho)$, if $\phi \in ker(d)$ (i.e. $d \phi =0$), then $\phi$ is a constant equivariant map.
\end{proposition} 

\begin{proof}
If $d \phi =0$, then by definition for every $(u,v) \in \Sigma (1,\Gamma)$ we have that $\phi (u)=\phi (v)$ and from equivariance $\phi (u)=\phi (v)$ for all $(u,v) \in \Sigma (1)$. Since $X$ is connected (because it is contractible) we get that for every $u,v \in \Sigma (0)$, $\phi (u) = \phi (v)$.
\end{proof}

\begin{remark}
Note that we used only the fact that $X$ is connected and not the fact that it is contractible, i.e. the proposition is true whenever $X$ is connected.
\end{remark}

\subsection{Localization}

Let $\left(v_{0},...,v_{j}\right)=\tau\in\Sigma(j)$, denote by $X_{\tau}$
the \emph{link} of $\tau$ in $X$, that is, the (pure) complex of dimension
$n-j-1$ consisting on simplices $\sigma=\left(w_{0},...,w_{k}\right)$
such that $\left\{ v_{0},...,v_{j}\right\} ,\left\{ w_{0},...,w_{k}\right\}$ are disjoint as sets
and $\left(v_{0},...,v_{j},w_{0},...,w_{k}\right)=\tau\sigma\in\Sigma(j+k+1)$.
The isotropy group $\Gamma_{\tau}$ acts by automorphisms on $X_{\tau}$
and if we denote by $\rho_{\tau}$ the restriction of $\rho$ to $\Gamma_{\tau}$,
we get that $\rho_{\tau}$ is a unitary representation of $\Gamma_{\tau}$.
In this section, we will always assume that $j \geq r$, which means that $X_\tau$ is a finite simplicial complex.
\begin{enumerate}
\item For $0\leq k\leq n-j-1$, denote by $\Sigma_{\tau}(k)$ the set of
ordered $k$-simplices and choose a set $\Sigma_{\tau}(k,\Gamma_{\tau})\subseteq\Sigma_\tau (k)$
of representatives of $\Gamma_{\tau}$-orbits.
\item For a simplex $\sigma\in\Sigma_{\tau}(k)$ denote by $m_{\tau}(\sigma) = m (\tau \sigma )$
and by our definition, $m_{\tau}(\sigma)> 0$ for every $\sigma$. Denote by $C_{\tau,k}$ the constant such that for every $\sigma \in  \Sigma_{\tau}(k)$ one has 
$$ \sum_{\gamma \in \Sigma_\tau (k+1), \sigma \subset \gamma} m_\tau (\gamma) = (k+2)! C_{\tau,k} m_\tau ( \sigma )$$
\item For a simplex $\sigma\in\Sigma_{\tau}(k)$, denote by $\Gamma_{\tau\sigma}$
the stabilizer of $\sigma$ in $\Gamma_{\tau}$ (or the stabilizer of $\tau\sigma$ in $\Gamma$). 
\item For $0\leq k\leq n-j-1$, denote by $C^{k}(X_{\tau},\rho_{\tau})$
the space of simplicial $k$-cochains of $X_{\tau}$ which are twisted
by $\rho_{\tau}$.
\item On $C^{k}(X_{\tau},\rho_{\tau})$ we define an Hermitian form as before, i.e. for every $\phi , \psi \in C^{k}(X_{\tau},\rho_{\tau})$:
$$
\left\langle \phi,\psi\right\rangle :=\sum_{\sigma\in\Sigma_\tau (k,\Gamma_\tau)}\frac{m_\tau (\sigma)}{(k+1)!\left|\Gamma_{\tau \sigma}\right|}\left\langle \phi(\sigma),\psi(\sigma)\right\rangle $$
note that for every $\phi \in C^{k}(X_{\tau},\rho_{\tau})$ we have:
$$ \Vert \phi \Vert^2 = \sum_{\sigma\in\Sigma_\tau (k,\Gamma_\tau)}\frac{m_\tau (\sigma)}{(k+1)!\left|\Gamma_{\tau \sigma}\right|}\vert \phi(\sigma) \vert^2 = \dfrac{1}{\vert \Gamma_\tau \vert} \sum_{\sigma\in\Sigma_\tau (k)}\frac{m (\tau \sigma)}{(k+1)!}\vert \phi(\sigma) \vert^2$$
\item On $C^{k}(X_{\tau},\rho_{\tau})$ the differential is defined as before and denoted by $d_\tau$ (and $\delta_\tau = d_\tau^*,\Delta^+_\tau = \delta_\tau d_\tau$). 

\end{enumerate}

\begin{proposition}
For $\tau$ of dimension $j \geq r$ for every $0\leq k\leq n-j-2$, we have $C_{\tau ,k} = C_{j+k+1}$
\end{proposition}

\begin{proof}
For $\sigma \in \Sigma_\tau (k)$ we have by definition
$$ (k+2)!C_{\tau,k} m_\tau (\sigma ) = \sum_{\begin{array}{c}
{\scriptstyle \gamma \in \Sigma_\tau (k+1)}\\
{\scriptstyle \sigma \subset \gamma }\end{array} } m_\tau (\gamma) = $$
$$ = \sum_{\begin{array}{c}
{\scriptstyle \gamma \in \Sigma_\tau (k+1)}\\
{\scriptstyle \sigma \subset \gamma }\end{array}} m (\tau \gamma) = \sum_{\begin{array}{c}
{\scriptstyle \eta \in \Sigma (j+k+2)}\\
{\scriptstyle \tau \sigma \subset \eta }\end{array}  } \dfrac{ (k+2)!}{(j+k+3)!} m (\eta) =$$
$$=  (k+2)! C_{j+k+1} m (\tau \sigma ) =  (k+2)! C_{j+k+1} m_\tau ( \sigma )$$ 

\end{proof}

Define the \emph{localization map} \[
C^{k}(X,\rho)\rightarrow C^{k-j-1}(X_{\tau},\rho_{\tau}),\;\phi\rightarrow\phi_{\tau}\]
Where $\phi_\tau$ is defined by $\phi_{\tau}(\sigma)=\phi(\tau\sigma)$.

\begin{proposition}
For every $\phi \in C^k (X, \rho )$, $k > r$ one has:
\begin{enumerate}
\item $$\sum_{\tau \in \Sigma(k-1, \Gamma)} \Vert \phi_\tau \Vert^2 = (k+1)! \Vert \phi \Vert^2$$
\item $$\sum_{\tau \in \Sigma(k-1, \Gamma)} \Vert \phi_\tau^0 \Vert^2 = \dfrac{k!}{C_{k-1}} \Vert \delta \phi \Vert^2$$
where $\phi_\tau^0$ is the projection of $\phi_\tau$ on the space of constant functions.
\end{enumerate}
\end{proposition}

\begin{proof}
\begin{enumerate}
\item This proposition is proven in \cite[Lemma 1.10]{BS} and since the proof doesn't take the weight function into account, we will not repeat it here.
\item For every $\phi \in C^k (X, \rho )$ and every $\tau \in \Sigma (k-1,\Gamma)$ one has:
$$ \phi_\tau^0 = \dfrac{1}{\sum_{v \in \Sigma_\tau (0)} m (\tau v)} \sum_{v \in \Sigma_\tau (0)} m (\tau v) \phi (\tau v) $$ 
Note that 
$$ \sum_{v \in \Sigma_\tau (0)} m (\tau v) = \dfrac{1}{(k+1)!} \sum_{\sigma \in \Sigma (k), \tau \subset \sigma} m (\sigma) = C_{k-1} m( \tau )$$
therefore we get 
$$ \phi_\tau^0  = \dfrac{1}{C_{k-1} m( \tau )} \sum_{v \in \Sigma_\tau (0)} m (\tau v) \phi (\tau v) = \dfrac{ (-1)^k}{C_{k-1}} \delta (\tau ) $$
and therefore
$$ \Vert \phi_\tau^0 \Vert^2 = \dfrac{1}{\vert \Gamma_\tau \vert } \sum_{v \in \Sigma_\tau (0)} m( \tau v) \vert \dfrac{(-1)^k}{C_{k-1}} \delta (\tau ) \vert^2 = $$
$$ = \dfrac{1}{\vert \Gamma_\tau \vert }   \dfrac{m( \tau )}{C_{k-1}} \vert \delta (\tau ) \vert^2 $$ 
and the equality in the proposition follows.
\end{enumerate}
 
\end{proof}

We shall also use the following proposition which is proven in \cite{BS} (the proposition is general and does not depend on the weight function, so we will not prove it here):

\begin{proposition}
\label{LaplacianLowerBound}
(\cite[Lemma 2.3]{BS}) For a simplex $\tau$ of dimension $j \geq r$, let $\alpha^{+}$ is the smallest positive eigenvalue of $\Delta_\tau^{+}$ on $C^{0}(X_\tau,\mathbb{R})$, which is the space of
(untwisted) $0$-cochains on $X_\tau$ with values in $\mathbb{R}$ (this can be thought of the space of representations in $\mathbb{R}$, where the group is the trivial group and not $\Gamma_\tau$). 
Then $\left\Vert \Delta_\tau^{+}\phi\right\Vert \geq\alpha^{+}\left\Vert \phi\right\Vert $
for all $\phi\in C^{0}(X_\tau,\rho_\tau)$ perpendicular to $ker(\Delta_\tau^{+})$ (Note that $X_\tau$ and $\Gamma_\tau$ are compact). \\
From now on denote by $\lambda (X_\tau)=\alpha^+$ the smallest positive eigenvalue of $\Delta_\tau^+$ on $C^{0}(X_\tau,\mathbb{R})$.

\end{proposition} 

\subsection{Restriction}
In this subsection we assume that $r=0$, i.e., we assume that $X$ is locally finite and that the action of $\Gamma$ is proper.

\begin{definition}
For $\phi \in C^k (X, \rho)$ and $\tau \in \Sigma (l)$ s.t. $k+l+1 \leq n$, the restriction of $\phi$ to $X_\tau$  is a function $\phi^\tau \in C^k (X_\tau, \rho_\tau)$ defined as follows: 
$$ \forall \sigma \in \Sigma_\tau (k), \phi^\tau (\sigma) = \phi (\sigma)$$
\end{definition}

\begin{proposition}
\label{restNorm1}
Let $\phi \in C^k (X,\rho)$ then
$$C_k \Vert \phi \Vert^2 = \sum_{u \in \Sigma (0, \Gamma )} \Vert \phi^u \Vert^2$$
In particular, if $X$ is finite and $\Gamma$ is trivial then
$$C_k \Vert \phi \Vert^2 = \sum_{u \in \Sigma (0)} \Vert \phi^u \Vert^2$$
\end{proposition}

\begin{proof}
$$ \sum_{u \in \Sigma (0, \Gamma)} \Vert \phi^u \Vert^2 = \sum_{u \in \Sigma (0, \Gamma)} \sum_{\tau \in \Sigma_u (k, \Gamma_u)} \dfrac{m_u (\tau)}{(k+1)! \vert \Gamma_{u \tau} \vert}\vert \phi^u (\tau) \vert^2 = $$
$$ = \sum_{u \in \Sigma (0, \Gamma)} \dfrac{1}{(k+1)! \vert \Gamma_u \vert }\sum_{\tau \in \Sigma_u (k)} m (u \tau)\vert \phi (\tau) \vert^2 = $$
$$ =  \sum_{u \in \Sigma (0, \Gamma)} \dfrac{1}{(k+1)! \vert \Gamma_u \vert }\sum_{\gamma \in \Sigma (k+1),u \subset \gamma} \dfrac{1}{k+2} m (\gamma)\vert \phi (\gamma - u) \vert^2 $$
where $\gamma - u$ means deleting the vertex of $u$ from $\gamma$. Changing the order of summation gives
$$ \sum_{\gamma \in \Sigma (k+1, \Gamma)} \dfrac{m(\gamma)}{(k+2)! \vert \Gamma_\gamma \vert }\sum_{u \in \Sigma (0),u \subset \gamma} \vert \phi (\gamma - u) \vert^2 = $$
$$ = \sum_{\gamma \in \Sigma (k+1, \Gamma)} \dfrac{m(\gamma)}{(k+2)! \vert \Gamma_\gamma \vert }\sum_{\sigma \in \Sigma (k),\sigma \subset \gamma} \dfrac{1}{(k+1)!}\vert \phi (\sigma) \vert^2 = $$
$$ = \sum_{\sigma \in \Sigma (k, \Gamma)} \dfrac{\vert \phi (\sigma) \vert^2}{(k+2)!(k+1)! \vert \Lambda_\sigma \vert }\sum_{\gamma \in \Sigma (k+1),\sigma \subset \gamma} m(\gamma)$$
so we get
$$ C_k \sum_{\sigma \in \Sigma (k, \Gamma)} \dfrac{m(\sigma) \vert \phi (\sigma) \vert^2}{(k+1)! \vert \Gamma_\sigma \vert } =  C_k \Vert \phi \Vert^2$$
\end{proof}

\begin{proposition}
\label{restNorm2}
Let $\phi \in C^0 (X,\rho)$ then
$$ C_1 \Vert d \phi \Vert^2 = \sum_{u \in \Sigma (0, \Gamma)} \Vert d_u \phi^u \Vert^2$$
where $d_u$ is the restriction of $d$ to the link of $u$. In particular, if $X$ is finite and $\Gamma$ is trivial we get that
$$ C_1 \Vert d \phi \Vert^2 = \sum_{u \in \Sigma (0)} \Vert d_u \phi^u \Vert^2$$
\end{proposition}

\begin{proof}
Note that 
$$ \forall (v_0,v_1) \in \Sigma_u (1), d_u \phi^u ((v_0,v_1)) = \phi (v_0)-\phi (v_1) = d \phi ((v_0,v_1)) = (d \phi)^u ((v_0,v_1)) $$
therefore $d_u (\phi^u) = (d \phi)^u$  and the proposition follows. 
\end{proof}

\begin{remark}
In this article, we use the above results only in the case where $X$ is finite, $\Gamma$ is trivial and $\rho$ is into $\mathbb{R}$. Since these assumptions doesn't significantly simplify the proofs we kept the more general setting above.
\end{remark}

\section{Criteria for vanishing cohomology}

This section has two parts - first we will give a geometrical criterion for vanishing of cohomology which generalizes the criteria given \cite{BS} to the framework of weighted simplicial complexes (see exact statement in corollary \ref{thmcriterion2}). Second, we will show that the fulfilment of the criterion of corollary \ref{thmcriterion2} for some $k \geq r+1$ implies that it is fulfilled for every $j$ s.t. $k \geq j \geq r+1$. \\ \\
For the convenience of the reader, we recall our setting and notations:
\begin{itemize}
\item  $X$ is a $n$-dimensional simplicial complex with the following properties:
\begin{enumerate}
\item $X$ is pure $n$-dimensional, i.e., every simplex in $X$ is contained in at least one $n$-dimensional simplex.
\item X is contractible.
\item all the links of $X$ of dimension $>0$ are connected.
\item There is a constant integer $ 0 \leq r \leq n-2$ such that all the links of simplices of dimension $r$ are finite.
\item $X$ is has a weight function $m$ - see definition \ref{weightedDef} above.
\end{enumerate}
\item $\Gamma$ is a locally compact, unimodular group of automorphisms of $X$ acting cocompactly on $X$ such that for every simplex $\sigma$ of dimension $r$ the stabilizer group of $\sigma$, denoted by $\Gamma_\sigma$, is compact. 
\item $\rho$ is a unitary representation of $\Gamma$ on a complex Hilbert space $H$.
\item For a finite link $X_\tau$, $\lambda (X_\tau)$ is the smallest positive eigenvalue of the untwisted Laplacian on $C^0 (X_\tau, \mathbb{R})$ - see proposition \ref{LaplacianLowerBound} above. \\ 
\end{itemize}

\subsection{A criterion for cohomology vanishing}

\begin{lemma}
Let $X$, $\Gamma$, $\rho$ be as above. For every $r+1 \leq k \leq n-1$ and every $\phi \in C^k (X,\rho)$ we have that
$$ k! \Vert d \phi \Vert^2 = k! \sum_{\tau \in \Sigma (k-1,\Gamma)} \left( \Vert d_\tau \phi_\tau \Vert^2 - C_k \dfrac{ k}{k+1} \Vert \phi_\tau \Vert^2 \right)$$
\end{lemma}

\begin{proof}
For $(v_0,...,v_{k+1}) = \sigma \in \Sigma (k+1)$ and $ 0 \leq i < j \leq k+1$ denote 
$$\sigma_{ij} = (v_0,...\hat{v_{i}},...,\hat{v_{j}},...,v_{k+1})$$ 
Then for every $\phi \in C^k (X ,\rho)$ we have
$$ \vert d \phi (\sigma) \vert^2 = \sum_{0 \leq i < j \leq k+1} \vert \phi_{\sigma_{ij}} (v_i) - \phi_{\sigma_{ij}} (v_j) \vert^2 - k \sum_{0 \leq i \leq k+1} \vert \phi (\sigma_i) \vert^2  = $$
$$ = \sum_{0 \leq i < j \leq k+1} \left( \vert \phi_{\sigma_{ij}} (v_i) - \phi_{\sigma_{ij}} (v_j) \vert^2 - \dfrac{k}{k+1} ( \vert \phi_{ \sigma_{ij}} (v_i) \vert^2+  \vert \phi_{ \sigma_{ij}} (v_j) \vert^2) \right) = $$
$$ =  \dfrac{1}{k!} \sum_{\tau \in \Sigma (k-1), \tau \subset \sigma} \left( \vert d_\tau \phi_{\tau} (\sigma - \tau) \vert^2 - \dfrac{k}{k+1} \sum_{v \in \Sigma_\tau (0), v \subset \sigma - \tau} \vert \phi_{\tau} (v) \vert^2 \right)$$
where $\sigma - \tau$ is the $1$-dimensional simplex obtained by deleting the the vertices of $\tau$ from $\sigma$. \\
We can use this equality to connect $\Vert d \phi \Vert$ to $\Vert d_\tau \phi_\tau \Vert$ and $\Vert \phi_\tau \Vert$:
$$ k! \Vert d \phi \Vert^2 = \sum_{\sigma \in \Sigma (k+1, \Gamma)} \dfrac{m( \sigma )}{(k+2)! \vert \Gamma_\sigma \vert} \vert d \phi (\sigma) \vert^2 = $$
$$  = \sum_{\sigma \in \Sigma (k+1, \Gamma)} \dfrac{1}{(k+2)! \vert \Gamma_\sigma \vert} \sum_{\tau \in \Sigma (k-1), \tau \subset \sigma} m_\tau ( \sigma - \tau ) \left( \vert d_\tau \phi_{\tau} (\sigma - \tau) \vert^2 - \dfrac{k}{k+1} \sum_{v \in \Sigma_\tau (0), v \subset \sigma - \tau} \vert \phi_{\tau} (v) \vert^2 \right) = $$
$$ = \sum_{\tau \in \Sigma (k-1, \Gamma)} \dfrac{1}{(k+2)! \vert \Gamma_\tau \vert} \sum_{\sigma \in \Sigma (k+1), \tau \subset \sigma} m_\tau ( \sigma - \tau ) \left( \vert d_\tau \phi_{\tau} (\sigma - \tau) \vert^2 - \dfrac{k}{k+1} \sum_{v \in \Sigma_\tau (0), v \subset \sigma - \tau} \vert \phi_{\tau} (v) \vert^2 \right) = $$
$$ = \sum_{\tau \in \Sigma (k-1, \Gamma)} \dfrac{1}{(k+2)! \vert \Gamma_\tau \vert} \dfrac{(k+2)!}{2!} \sum_{\eta \in \Sigma_\tau (1)} m_\tau ( \eta ) \left( \vert d_\tau \phi_{\tau} (\eta) \vert^2 - \dfrac{k}{k+1} \sum_{v \in \Sigma_\tau (0), v \subset \eta} \vert \phi_{\tau} (v) \vert^2 \right)=$$
$$ =  \sum_{\tau \in \Sigma (k-1, \Gamma)} \dfrac{1}{\vert \Gamma_\tau \vert}  \sum_{\eta \in \Sigma_\tau (1)} \dfrac{m_\tau ( \eta ) }{2}\vert d_\tau \phi_{\tau} (\eta) \vert^2 - $$
$$ - \dfrac{k}{k+1} \sum_{\tau \in \Sigma (k-1, \Gamma)} \dfrac{1}{\vert \Gamma_\tau \vert} \sum_{\eta \in \Sigma_\tau (1)}  \sum_{v \in \Sigma_\tau (0), v \subset \eta} \dfrac{m_\tau (\eta )}{2}\vert \phi_{\tau} (v) \vert^2$$
Note that 
$$ \sum_{\tau \in \Sigma (k-1, \Gamma)} \dfrac{1}{\vert \Gamma_\tau \vert}  \sum_{\eta \in \Sigma_\tau (1)} \dfrac{m_\tau ( \eta ) }{2}\vert d_\tau \phi_{\tau} (\eta) \vert^2 =  \sum_{\tau \in \Sigma (k-1, \Gamma)} \Vert d_\tau \phi_\tau \Vert^2$$
and that 

$$ \dfrac{k}{k+1} \sum_{\tau \in \Sigma (k-1, \Gamma)} \dfrac{1}{\vert \Gamma_\tau \vert} \sum_{\eta \in \Sigma_\tau (1)}  \sum_{v \in \Sigma_\tau (0), v \subset \eta} \dfrac{m_\tau (\eta )}{2}\vert \phi_{\tau} (v) \vert^2 = $$
$$=  \dfrac{k}{k+1} \sum_{\tau \in \Sigma (k-1, \Gamma)} \dfrac{1}{\vert \Gamma_\tau \vert} \sum_{v \in \Sigma_\tau (0)} \vert \phi_{\tau} (v) \vert^2  \sum_{\eta \in \Sigma_\tau (1), v \subset \eta} \dfrac{m_\tau (\eta )}{2} =$$
$$ = \dfrac{k}{k+1} \sum_{\tau \in \Sigma (k-1, \Gamma)} \dfrac{1}{\vert \Gamma_\tau \vert} \sum_{v \in \Sigma_\tau (0)} \vert \phi_{\tau} (v) \vert^2 m_\tau (v) C_{\tau, 0}  =  \sum_{\tau \in \Sigma (k-1, \Gamma)} C_k \dfrac{k}{k+1} \Vert \phi_\tau \Vert^2  $$
So we get the desired equality.

\end{proof}

\begin{theorem}
\label{thmcriterion}
Let $X$, $\Gamma$ and $\rho$ be as above and let $r+1 \leq k \leq n-1$. Denote $\lambda_k = \min_{\tau \in  \Sigma (k-1)} \lambda (X_\tau)$ and assume that $\lambda_k > C_k \frac{k}{k+1}$. Then there is an $\varepsilon > 0$ such that for every $\phi \in C^k (X ,\rho), d \phi =0$ the following holds: 
$$\Vert \delta \phi \Vert^2 \geq \varepsilon \Vert \phi \Vert^2$$

\end{theorem}     

\begin{proof}
By the above lemma, for every $\phi \in C^k (X ,\rho)$ we have
$$ k! \Vert d \phi \Vert^2 = \sum_{\tau \in \Sigma (k-1,\Gamma)} \left( \Vert d_\tau \phi_\tau \Vert^2 - C_k \dfrac{ k}{k+1} \Vert \phi_\tau \Vert^2 \right)$$
For every $\tau \in \Sigma (k-1)$ link $X_\tau$ is connected and finite so the kernel of $d_\tau$ consists of constant maps and we get that
$$ \Vert d_\tau \phi_\tau \Vert^2 \geq \lambda_k \Vert \phi_\tau \Vert^2 - \lambda_k \Vert \phi_\tau^0 \Vert^2$$
Therefore
$$ k! \Vert d \phi \Vert^2 \geq \sum_{\tau \in \Sigma (k-1,\Gamma)} \left( (\lambda_k - C_k \dfrac{ k}{k+1}) \Vert \phi_\tau \Vert^2 - \lambda_k \Vert \phi_\tau^0 \Vert^2 \right) $$
Recall that
$$\sum_{\tau \in \Sigma(k-1, \Gamma)} \Vert \phi_\tau \Vert^2 = (k+1)! \Vert \phi \Vert^2 $$
and that 
$$\sum_{\tau \in \Sigma(k-1, \Gamma)} \Vert \phi_\tau^0 \Vert^2 = \dfrac{k!}{C_{k-1}} \Vert \delta \phi \Vert^2$$
So we get that
$$ k! \Vert d \phi \Vert^2 \geq (k+1)! (\lambda_k - C_k \dfrac{ k}{k+1}) \Vert \phi \Vert^2 - \lambda_k \dfrac{k!}{C_{k-1}} \Vert \delta \phi \Vert^2 $$  
For $\phi \in ker (d)$ this reads
$$ \Vert \delta \phi \Vert^2  \geq \dfrac{C_{k-1}(k+1)}{\lambda_k} (\lambda_k - C_k \dfrac{ k}{k+1}) \Vert \phi \Vert^2$$
So for $\varepsilon = \dfrac{C_{k-1}(k+1)}{\lambda_k} (\lambda_k - C_k \dfrac{ k}{k+1}) > 0$ we get the desired inequality.

\end{proof}

\begin{remark}
One should note that in the proof of the above theorem (and the preceding lemma) the fact that $X$ is contractible was not used and that this theorem holds even if $X$ is not contractible. 
\end{remark}

As corollary of the above theorem and theorem \ref{CohomologyVanish} we get the following corollary which is a version of \cite[Theorem 2.5]{BS} in our setting:

\begin{corollary}
\label{thmcriterion2}
Let $X$, $\Gamma$ and $\rho$ be as above and let $r+1 \leq k \leq n-1$. Denote $\lambda_k = \min_{\tau \in  \Sigma (k-1)} \lambda (X_\tau)$ and assume that $\lambda_k > C_k \frac{k}{k+1}$. Then $H^k (X, \rho) = 0$. Also, if $r=0$ (i.e. $X$ is locally finite and the action of $\Gamma$ is proper), then $H^k (X, \rho)= H^k (\Gamma, \rho)= 0$.
\end{corollary}

\subsection{Laplacian eigenvalues estimations}

Recall that for a finite simplicial complex $Y$ we denoted by $\lambda (Y)$ the smallest positive eigenvalue of the untwisted Laplacian on $Y$ (see proposition \ref{LaplacianLowerBound}). In this section we will show the following theorem about (untwisted) Laplacian eigenvalues (note that we don't use the group action in any way and this is purely a combinatorial result):

\begin{theorem}
\label{thmeigenval}
Let $X$ be as above and let $r+1 \leq k \leq n-1$. Denote $\lambda_k = \inf_{\tau \in \Sigma (k-1)} \lambda (X_\tau)$ and assume that $\lambda_k > C_k \frac{k}{k+1}$, then for every $j$ such that $r+1 \leq j < k$ we have that $ \inf_{\tau \in \Sigma (j-1)} \lambda (X_\tau) = \lambda_j > C_j \frac{j}{j+1}$.
\end{theorem}

Combined with corollary \ref{thmcriterion2} this implies the following corollary which is the first implication of theorem \ref{introthm} in the introduction when $m$ is the weight function given in \cite{BS} (see example \ref{BSweights} to see that $C_k = n-k$ in \cite{BS}):

\begin{corollary}
\label{cohomologyVanishResult}
Let $X$, $\Gamma$, $\rho$ be as above and let $r+1 \leq k \leq n-1$. Assume $\min_{\sigma \in \Sigma (k-1)} \lambda (X_\sigma)=\lambda_k > C_k \frac{k}{k+1}$, then for every $r+1 \leq j \leq k \leq n-1$ we have that $H^j (X, \rho)  = 0$. Also, if $r=0$ (i.e. $X$ is locally finite) then $\min_{\sigma \in \Sigma (k-1)} \lambda (X_\sigma)=\lambda_k > C_k \frac{k}{k+1}$ implies $H^j (\Gamma, \rho)  =H^j (X, \rho)= 0$ for every $1 \leq j \leq k$.
\end{corollary}

\begin{proposition}

\label{thmeval}
Let $X$ be as above, fix $r < j \leq n-2$ and $\tau \in \Sigma (j-1)$. If for every $v \in \Sigma_\tau (0)$, we have $\lambda( X_{\tau v}) \geq \mu$, then
$$ \lambda (X_\tau )  \geq  C_j \left( 2 - \dfrac{C_{j+1}}{\mu} \right)$$  
\end{proposition}  

\begin{proof}

Recall that for $\phi \in C^0 (X_\tau)$
$$\Delta_\tau^+ \phi (u) = C_{\tau,0} \phi (u) - \dfrac{1}{m_\tau (u)} \sum_{(u,v) \in \Sigma_\tau (1)} m_\tau((u,v)) \phi (v)$$

Take $\phi \in  C^0 (X_\tau, \mathbb{R})$ to be an eigenfunction of the untwisted Laplacian $\Delta^+_\tau$ with the eigenvalue $\lambda (X_\tau)$, i.e. 
$$\Delta_\tau^+ \phi (u) = \lambda (X_\tau) \phi (u)$$ 

Fix $v \in \Sigma_\tau (0)$ and denote by $( \phi^v)^0$ the projection of $\phi^v$ to the space of constant maps on $X_{\tau v}$ and by $( \phi^v)^1$ its orthogonal compliment. Explicitly we have
$$ \sum_{u \in \Sigma_{\tau v} (0)} m_{\tau v} (u) \phi^v (u) = \sum_{u \in \Sigma_{\tau v} (0)} m_{\tau v} (u) (\phi^v)^0 = (\phi^v)^0 \sum_{u \in \Sigma_{\tau v} (0)} m_{\tau v} (u)= $$
$$   =(\phi^v)^0 \sum_{u \in \Sigma_\tau (0), (v,u) \in \Sigma_\tau (1)} m_\tau ((v,u)) =(\phi^v)^0 C_{\tau , 0} m_\tau (v)$$
Therefore
$$ (\phi^v)^0 =\dfrac{1}{C_{\tau , 0} m_\tau (v)} \sum_{u \in \Sigma_{\tau v} (0)} m_{\tau v} (u) \phi^v (u)$$
and
$$ (\phi^v)^1 = \phi^v - (\phi^v)^0$$
Note that since $\Delta_\tau^+ \phi (v) = \lambda (X_\tau) \phi (v)$, we get that
$$ \lambda (X_\tau) \phi (v) = C_{\tau, 0} \phi (v) - \dfrac{1}{m_\tau (v)} \sum_{(v,u) \in \Sigma_\tau (1)} m_\tau ((v,u)) \phi (u) =$$
$$= C_{\tau, 0} \phi (v) - \dfrac{1}{m_\tau (v)} \sum_{u \in \Sigma_{\tau v} (0)} m_{\tau v} (u) \phi^v (u)$$
Therefore
$$ (\phi^v)^0 = \dfrac{C_{\tau,0} - \lambda (X_\tau)}{C_{\tau,0}} \phi (v)$$

Since $X_{\tau v}$ is connected for every $v \in \Sigma_\tau (0)$, the kernel of $\Delta^+_{\tau v}$ is the space of constant maps (see proposition \ref{ConstantEqui} and the remark that follows it). Therefore, by the definition of $\lambda (X_{\tau v})$, we have for every  $v \in \Sigma_\tau (0)$
$$ \Vert d_{\tau v} \phi^v \Vert^2 \geq \lambda (X_{\tau v}) \Vert (\phi^v)^1 \Vert^2 = \lambda (X_{\tau v}) \Vert (\phi^v) \Vert^2 - \lambda (X_{\tau v}) \Vert (\phi^v)^0 \Vert^2$$
Since $X_\tau$ is finite we can apply proposition \ref{restNorm2} in the untwisted case:
$$C_{\tau,1} \Vert d_\tau \phi \Vert^2 = \sum_{v \in \Sigma_\tau (0)} \Vert d_{\tau v} \phi^v \Vert^2$$
Combined with the above inequality this yields:
$$ C_{\tau,1} \Vert d_\tau \phi \Vert^2 \geq \sum_{v \in \Sigma_\tau (0)} \lambda (X_{\tau v}) \Vert (\phi^v) \Vert^2 - \lambda (X_{\tau v}) \Vert (\phi^v)^0 \Vert^2 $$
Note that 
$$\Vert (\phi^v)^0 \Vert^2 =  \sum_{u \in \Sigma_{\tau v} (0)} m_{\tau v} (u) \vert (\phi^v)^0 \vert^2 =$$
$$=   C_{\tau,0} m_\tau (v) \vert (\phi^v)^0 \vert^2 = m_\tau (v) \dfrac{ (C_{\tau,0} - \lambda (X_\tau))^2}{C_{\tau,0}} \vert \phi (v) \vert^2 $$  
So we get
$$ C_{\tau,1} \lambda (X_\tau) \Vert \phi \Vert^2 = C_{\tau,1} \Vert d \phi \Vert^2 \geq $$
$$ \geq \sum_{v \in \Sigma_\tau (0)} \lambda (X_{\tau v}) \left( \Vert (\phi^v) \Vert^2 -  m_\tau (v) \dfrac{ (C_{\tau,0} - \lambda (X_\tau))^2}{C_{\tau,0}} \vert \phi (v) \vert^2 \right) \geq $$
$$ \geq \sum_{v \in \Sigma_\tau (0)} \mu  \left( \Vert (\phi^v) \Vert^2 -  m_\tau (v) \dfrac{ (C_{\tau,0} - \lambda (X_\tau))^2}{C_{\tau,0}} \vert \phi (v) \vert^2 \right) =$$
$$=  \mu \left( C_{\tau ,0} - \dfrac{ (C_{\tau,0} - \lambda (X_\tau))^2}{C_{\tau,0}} \right) \Vert \phi \Vert^2$$
where the last equality is due to proposition \ref{restNorm1} in the untwisted case. \\
Therefore
$$ C_{\tau,1} \lambda (X_\tau) \geq   \mu  \left( C_{\tau ,0} - \dfrac{ (C_{\tau,0} - \lambda (X_\tau))^2}{C_{\tau,0}} \right)$$
Which yields
$$ \lambda (X_\tau) \left( \lambda (X_\tau) \dfrac{\mu}{C_{\tau,0}} - 2 \mu + C_{\tau,1} \right) \geq 0$$
and since $\lambda (X_\tau) > 0$ we get
$$ \lambda (X_\tau)\geq \dfrac{C_{\tau,0}(-C_{\tau,1} + 2 \mu)}{\mu}$$
Recall that $C_{\tau ,0} =  C_{j}$, $C_{\tau ,1} =  C_{j+1}$ and therefore
$$ \lambda (X_\tau) \geq \dfrac{ C_{j}(- C_{j+1} + 2 \mu)}{\mu} =  C_j \left( 2 - \dfrac{C_{j+1}}{\mu} \right)$$

\end{proof}

Now we can prove theorem \ref{thmeigenval}:

\begin{proof}
Let $r+1 < k \leq n-1$. Assume 
$$\inf_{\tau \in \Sigma(k-1)} \lambda (X_\tau) = \lambda_k > C_k \dfrac{k}{k+1}$$ 
We shall prove by induction that $\inf_{\tau \in \Sigma(j-1)} \lambda (X_\tau) = \lambda_j > C_j \frac{j}{j+1}$ (as long as $r < j \leq k$). \\
Fix $\tau \in \Sigma (j-1)$ and assume that that the theorem is true for $j+1$, i.e. that 
$$\inf_{\eta \in \Sigma(j)} \lambda (X_\eta) = \lambda_{j+1} > C_{j+1} \dfrac{j+1}{j+2}$$ 
From the definition of $\lambda_{j+1}$ it is clear that for every $v \in \Sigma_\tau (0)$, we have $\lambda (X_{\tau v}) \geq \lambda_{j+1}$. Therefore by proposition \ref{thmeval} 
$$\lambda (X_\tau) \geq C_j \left( 2 - \dfrac{C_{j+1}}{\lambda_{j+1}} \right) $$
and since this is true for every $\tau \in \Sigma (j-1)$, we have that
$$\lambda_j = \inf_{\tau \in \Sigma (j-1)} \lambda (X_\tau) \geq C_j \left( 2 - \dfrac{C_{j+1}}{\lambda_{j+1}} \right)> C_j \left( 2 - \dfrac{C_{j+1}}{C_{j+1} \frac{j+1}{j+2}} \right)= C_j \dfrac{j}{j+1}$$
where the sharp inequality is due to the induction assumption.
\end{proof}

\section{Property (T)}
In this section we shall show that that every one of the above criteria given for the vanishing of some cohomology in corollary \ref{thmcriterion2} implies property (T) even if $X$ is not locally finite - see exact formulation in theorem \ref{propertyTresult} below. The structure of this section is as follows: first, we shall prove a general criterion for property (T) in the case where the group acts on a $2$-dimensional simplicial complex that is not necessarily contractible - see exact formulation in theorem \ref{Prop(T)reduction}. Second, for $\Gamma$ acting on $X$ as above, we shall construct such a complex $X'$ using the group action on $X$. Then we'll show that if the criterion for the cohomology vanishing stated in corollary \ref{thmcriterion2} is met, then the criterion for property (T) is met.

\subsection{The $2$-dimensional, locally finite case}

In this subsection we shall prove a criterion for property (T) for groups acting on a $2$-dimensional simplicial complex (see exact formulation below). This criterion is well known (see for instance \cite{BHV}[Theorem 5.5.2] or \cite{Gromov}[Subsection 3.11]) and the proof is given here for completeness. \\
Let us recall the definition of property (T) (the reader not familiar with property (T) can find a much more complete discussion in \cite{BHV}):
\begin{definition}
A topological group $\Gamma$ is said to have Kazhdan property (T) if the following holds: there exists a constant $\kappa >0$ and a compact subset $Q \subset \Gamma$ such that for every unitary representation $\rho: \Gamma \rightarrow \mathcal{U} (H)$ (where $H$ is a complex Hilbert space) one of the following holds:
\begin{enumerate}
\item $\rho$ has a non trivial invariant vector.
\item For every unit vector $\xi \in H$ we have 
$$\max_{\gamma \in Q} \vert \rho (\gamma)\xi - \xi \vert \geq \kappa$$
\end{enumerate}
In such a case, the couple $(Q,\kappa)$ is called a Kazhdan pair of $\Gamma$.
\end{definition}

The aim of this subsection is to show a criterion for property (T) in the following setting: 
let $X'$ a $2$-dimensional simplicial complex with the following properties:
\begin{enumerate}
\item $X'$ is pure $2$-dimensional.
\item $X'$ is connected.
\item $X'$ is locally finite.
\item The $1$-dimensional links of $X'$ are connected.
\item $X'$ has a weight function $m'$.
\end{enumerate}
Let $\Gamma$ be a locally compact, unimodular group of automorphisms of $X'$ acting cocompactly on $X'$ such that the action is proper. \\
 In order to distinguish $X'$ from $X$, we shall denote the ordered $k$-simplices of $X'$ as $\Sigma ' (k)$ and the orbit representatives as $\Sigma ' (k,\Gamma)$. Also, the constants of $m'$ will be denoted as $C'_k$. 
\begin{remark}
The reader should note that the assumptions on $X'$ are different from the assumptions on $X$. For instance, $X'$ is assumed to be locally finite and this is not the case in $X$ if $r>0$. Also, $X$ is assumed to be contractible where $X'$ is only assumed to be connected.
\end{remark}

\begin{theorem}
\label{Prop(T)reduction}
Let $X'$ and $\Gamma$ as above.  If $\min_{u \in \Sigma '(0)} \lambda (X_u') = \lambda_1' >  \frac{C'_1}{2}$ then $\Gamma$ has property (T). 
\end{theorem}

In order to prove the above theorem we shall need the following lemma:
\begin{lemma}
Let $X'$ and $\Gamma$ as above. If there is an $\varepsilon >0$ such that for every unitary representation $\rho$ of $\Gamma$ without a non trivial invariant vector we have
$$\forall \psi \in C^0 (X', \rho), \Vert d \psi \Vert^2 \geq \varepsilon \Vert \psi \Vert^2$$
then $\Gamma$ has property (T). 
\end{lemma}

\begin{proof}
For every $v \in \Sigma ' (0) $ choose (only) one group element $s_v \in \Gamma$ such that $s_v  v \in \Sigma ' (0,\Gamma)$. Denote the vertices in $\Sigma ' (1,\Gamma)$ as
$$  V(\Sigma ' (1,\Gamma))= \lbrace v  \in \Sigma ' (0) : \exists u \in \Sigma ' (0), (u,v) \in \Sigma ' (1,\Gamma) \text{ or } (v,u) \in \Sigma ' (1,\Gamma) \rbrace$$
 and define $S \subset \Gamma$ as
$$S = \lbrace s_v \in \Gamma : v \in  V(\Sigma ' (1,\Gamma))   \rbrace$$
Note that $S$ is finite, since the group action is cocompact .Define $Q = \bigcup_{u \in \Sigma ' (0, \Gamma)} \Gamma_u \cup S \subset \Gamma$. $Q$ is compact as a finite union of compact sets. Denote
 $$\kappa = \dfrac{\sqrt{\varepsilon}}{\sqrt{\varepsilon} + \sqrt{8C_0'}}$$
We claim that $(Q,\kappa)$ is a Kazhdan pair of $\Gamma$. Let $\rho$ be a unitary representation of $\Gamma$ without a non trivial invariant vector. Let $\xi \in H$ be a unit vector, if 
$$\max_{\gamma \in Q} \vert \rho (\gamma) \xi - \xi \vert \geq 1$$
we are done because $\kappa <1$. Assume now that 
$$\max_{\gamma \in Q} \vert \rho (\gamma) \xi - \xi \vert < 1$$
For every $u \in \Sigma '(0,\Gamma)$ consider the closure of the convex hull of the orbit of $\xi$: $\overline{conv(\Gamma_u \xi)}$. This is a convex closed set of $H$ and therefore it has a unique point denoted by $\xi_u \in \overline{conv(\Gamma_u  \xi)}$ such that $\vert \xi_u \vert$ is minimal (this is a classical fact about convex closed sets in Hilbert spaces - see \cite{Rudin}[Theorem 4.10]). Since for every $\gamma \in \Gamma_u$ we have $\gamma  (\Gamma_u  \xi )= \Gamma_u  \xi$ we get that $\gamma (\overline{conv(\Gamma_u  \xi) })= \overline{conv(\Gamma_u  \xi)}$. From the uniqueness of $\xi_u$ we have that
 $$\forall \gamma \in \Gamma_u, \gamma \xi_u = \xi_u$$
 In other words, $\xi_u$ is fixed by $\Gamma_u$. Also, 
 $$\forall \xi ' \in \Gamma_u  \xi, \vert \xi ' - \xi \vert \leq \max_{\gamma \in Q} \vert \rho (\gamma) \xi - \xi \vert < 1$$
Since the norm of $H$ is convex we get that  
$$\forall \xi ' \in \overline{conv(\Gamma_u  \xi)}, \vert \xi ' - \xi \vert \leq \max_{\gamma \in Q} \vert \rho (\gamma) \xi - \xi \vert < 1$$
In particular, 
$$\vert \xi_u - \xi \vert \leq \max_{\gamma \in Q} \vert \rho (\gamma) \xi - \xi \vert  <1$$ 
and as a result 
$$\vert \xi_u \vert \geq 1-  \max_{\gamma \in Q} \vert \rho (\gamma) \xi - \xi \vert > 0$$ 
Define $\psi \in C^0 (X', \rho)$ by 
$$\forall v \in \Sigma ' (0), \psi (v) = \rho(s_v^{-1} ) \xi_{s_v v}$$
$\psi$ is equivariant: by the definition of $\Sigma ' (0, \Gamma)$ we have
$$ \forall \gamma \in \Gamma, \forall v \in \Sigma ' (0), s_{\gamma v} \gamma v = s_v v$$
and therefore for every $\gamma \in \Gamma$ and $v \in \Sigma ' (0)$ we have
$$ \rho (\gamma ) \psi (v ) = \rho (\gamma ) \rho(s_v^{-1} ) \xi_{s_v v} = \rho (s_{\gamma v}^{-1}) \rho (s_{\gamma v} \gamma s_v^{-1} ) \xi_{ s_{\gamma v} \gamma v} = \rho (s_{\gamma v}^{-1})  \xi_{ s_{\gamma v} \gamma v} = \psi (\gamma v)$$
where the third equality is due to the fact that $s_{\gamma v} \gamma s_v^{-1} \in \Gamma_{s_{\gamma v} \gamma v}$ and that $\Gamma_{s_{\gamma v} \gamma v}$ fixes $\xi_{ s_{\gamma v} \gamma v}$. \\
For every $v \in  V(\Sigma ' (1,\Gamma))$the following holds:
$$\vert \xi - \psi (v) \vert = \vert \xi - \rho( s_v^{-1}) \xi_{s_v v} \vert = \vert \rho( s_v) \xi - \xi_{s_v v} \vert \leq$$
$$\leq \vert \rho( s_v) \xi - \xi \vert + \vert \xi - \xi_{s_v v} \vert \leq 2 \max_{\gamma \in Q} \vert \rho (\gamma) \xi - \xi \vert$$
So for all $\sigma = (u,v) \in \Sigma ' (1, \Gamma)$ we have 
$$\vert d \psi (\sigma) \vert = \vert \psi (u) - \psi (v) \vert \leq  \vert \psi (u) - \xi \vert  + \vert \xi - \psi (v) \vert \leq 4 \max_{\gamma \in Q} \vert \rho (\gamma) \xi - \xi \vert$$
Now we can use the assumption $\Vert d \psi \Vert^2 \geq \varepsilon \Vert \psi \Vert^2$. First note that 
$$ \Vert \psi \Vert^2 = \sum_{u \in \Sigma ' (0,\Gamma)} \frac{m' (u)}{\left|\Gamma_{u}\right|} \vert \xi_u \vert^2 \geq \sum_{u \in \Sigma ' (0,\Gamma)} \frac{m' (u)}{\left|\Gamma_{u}\right|}  \left( 1-  \max_{\gamma \in Q} \vert \rho (\gamma) \xi - \xi \vert \right)^2$$
Also note that 
$$ \Vert d \psi \Vert^2 = \sum_{\sigma \in \Sigma ' (1, \Gamma)} \frac{m' (\sigma )}{2 \left|\Gamma_{\sigma}\right|}  \vert d \psi (\sigma ) \vert^2 \leq \left( 4 \max_{\gamma \in Q} \vert \rho (\gamma) \xi - \xi \vert \right)^2 \sum_{\sigma \in \Sigma ' (1, \Gamma)} \frac{m' (\sigma)}{2 \left|\Gamma_{\sigma}\right|} = $$
$$= \left( 4 \max_{\gamma \in Q} \vert \rho (\gamma) \xi - \xi \vert \right)^2 \sum_{\sigma \in \Sigma ' (1, \Gamma)} \frac{1}{ \left|\Gamma_{\sigma}\right|}  \sum_{u \in \Sigma ' (0), u \subset \sigma} \dfrac{m' (\sigma)}{4}  = $$
$$= \left( 4 \max_{\gamma \in Q} \vert \rho (\gamma) \xi - \xi \vert \right)^2 \sum_{u \in \Sigma ' (0, \Gamma)} \frac{1}{ \left|\Gamma_{u}\right|}  \sum_{\sigma \in \Sigma ' (1), u \subset \sigma} \dfrac{m' (\sigma)}{4}  = $$
$$= \left( 4 \max_{\gamma \in Q} \vert \rho (\gamma) \xi - \xi \vert \right)^2 \sum_{u \in \Sigma ' (0, \Gamma)} \frac{C_0' m' (u) }{2 \left|\Gamma_{u}\right|} $$
By the inequality $\Vert d \psi \Vert^2 \geq \varepsilon \Vert \psi \Vert^2$ we get
$$\left( 4 \max_{\gamma \in Q} \vert \rho (\gamma) \xi - \xi \vert \right)^2 \dfrac{C_0'}{2} \sum_{u \in \Sigma ' (0, \Gamma)} \frac{ m' (u) }{ \left|\Gamma_{u}\right|}  \geq \varepsilon \left( 1-  \max_{\gamma \in Q} \vert \rho (\gamma) \xi - \xi \vert \right)^2  \sum_{u \in \Sigma ' (0,\Gamma)} \frac{m' (u)}{\left|\Gamma_{u}\right|} $$ 
which yields
$$ \max_{\gamma \in Q} \vert \rho (\gamma) \xi - \xi \vert \geq \dfrac{\sqrt{\varepsilon}}{\sqrt{\varepsilon} + \sqrt{8C_0'}} = \kappa$$
so $\left( Q, \kappa \right)$ is a Kazhdan pair.
\end{proof}

Now we can prove theorem \ref{Prop(T)reduction}:
\begin{proof}
Let $X'$, $\Gamma$ as above and let $\rho$ be a unitary representation of $\Gamma$ without a non trivial invariant vector. \\
To distinguish the differentials in this proof we shall denote 
$$d_0 : C^0 (X' ,\rho  ) \rightarrow C^1 (X' ,\rho  )$$
$$ d_1 : C^1 (X' ,\rho ) \rightarrow C^2 (X' ,\rho  )$$
and $\delta = d_1^*$. \\
Denote $\lambda_1' = \min_{u \in \Sigma ' (0)} \lambda (X_u')$ and assume that $\lambda_1' >  \frac{C'_1}{2}$. By theorem \ref{thmcriterion}, there is an $\varepsilon > 0$ such that for every $\phi \in C^1 (X ,\rho), d_1 \phi =0$, the following inequality holds
$$\Vert \delta \phi \Vert^2 \geq \varepsilon \Vert \phi \Vert^2$$ 
For every $\psi \in C^0 (X', \rho)$ we can take $d_0 \psi = \phi \in ker(d_1)$. Since $\delta d_0 = \Delta^+$ we get
$$ \Vert \Delta^+ \psi \Vert^2 \geq \varepsilon \langle \Delta^+ \psi, \psi \rangle$$
or
$$ \langle d_0 (\Delta^+)^{\frac{1}{2}} \psi, d_0 (\Delta^+)^{\frac{1}{2}} \psi \rangle \geq \varepsilon \langle (\Delta^+)^{\frac{1}{2}} \psi,(\Delta^+)^{\frac{1}{2}} \psi \rangle$$
Note that $ker ((\Delta^+)^{\frac{1}{2}})= \lbrace 0 \rbrace$, because $ker (\Delta^+)^{\frac{1}{2}} = ker (\Delta^+) = ker (d_0)$, and since $X'$ is connected $ker (d_0)$ contains only constant equivariant maps (see proposition \ref{ConstantEqui} and the remark that follows it), but $\rho$ has no non trivial invariant vector and therefore  $ker (d_0) = \lbrace 0 \rbrace$.  So we get 
$$\overline{im ((\Delta^+)^{\frac{1}{2}})} = (ker ((\Delta^+)^{\frac{1}{2}}))^\perp=\lbrace 0 \rbrace^\perp = C^{0}(X',\rho)$$
and therefore for every $\psi \in  C^{0}(X',\rho)$ the inequality
$$ \Vert d_0 \psi \Vert^2 \geq \varepsilon \Vert \psi \Vert^2$$
holds and we are done by the above lemma.
\end{proof}

As we shall see in the next subsection, it is sometimes useful to state the Laplacian eigenvalue condition in theorem \ref{Prop(T)reduction} in matrix form. To do this we'll need the following lemma (which is a standard result concerning normalized Laplacian adapted to the weighted Laplacian that we use):

\begin{lemma}
\label{MatrixReduction}
Let $Z=(V,E)$ be a connected finite graph with a weight function $m$. And $\Delta^+: C^0 (Z, \mathbb{R}) \rightarrow C^0 (Z, \mathbb{R})$ is defined as before (with respect to the weight function). Define the following $\vert V \vert \times \vert V \vert$ matrix:
$$ (A_Z) (u,v) = \begin{cases} C_0 & u= v \\
								\dfrac{-m((u,v))}{\sqrt{m(u)m(v)}} & (u,v) \in E \\
								0 & \text{otherwise}
					\end{cases} $$     
Then the eigenvalues of $\Delta^+$ and of the matrix $A_Z$ coincide and in particular $\lambda (Z)$ is equal to the smallest positive eigenvalue of $A_Z$. 		
\end{lemma}

\begin{proof}
Let $x=(x_v)$ be an eigenvector of $A_Z$ with eigenvalue $\mu$. Define a new vector $y=(y_v)$ as $y_v = \sqrt{m(v)} x_v$, then $A_Z x = \mu x$ is equivalent to 
$$ \forall v \in V, C_0 x_v - \sum_{u \in V, (v,u) \in E} \dfrac{m((u,v))}{\sqrt{m(u)m(v)}} x_u = \mu x_v$$
which is equivalent to
$$ \forall v \in V, C_0 \dfrac{1}{\sqrt{m(v)}} y_v - \sum_{u \in V, (v,u) \in E} \dfrac{m((u,v))}{m(u) \sqrt{m(v)}} y_u = \mu \dfrac{1}{\sqrt{m(v)}} y_v$$
Multiplying by $\sqrt{m(v)}$ gives $ \Delta^+ y = \mu y$ (and in the same way we can start with an eigenvector $y$ of $\Delta^+$ and get an eigenvector of $A_Z$).

\end{proof}

\subsection{Criterion for property (T) in the non locally finite case} 

In this subsection we shall prove that if the criterion for the vanishing of $H^k (X,\rho)$ given in corollary \ref{thmcriterion2} holds, then $\Gamma$ has property (T). Specifically:

\begin{theorem}
\label{propertyTresult}
Let $X$, $\Gamma$ as in section 3 and let $r+1 \leq k \leq n-1$. If $\lambda_k = \min_{\tau \in \Sigma (k-1)} \lambda (X_\tau) > C_k \dfrac{k}{k+1}$, then $\Gamma$ has property (T).
\end{theorem}

\begin{remark}
Note that the theorem stated above implies the the second implication of Theorem \ref{introthm} in the introduction when $m$ is taken to be the weight function given in \cite{BS} (see example \ref{BSweights} to see that with the weight function used in \cite{BS} we have $C_k = n-k$).
\end{remark}

\begin{remark}
If the link of every vertex in $X$ is compact and the action of $\Gamma$ is proper, this theorem follows directly from corollary \ref{cohomologyVanishResult}.
\end{remark}

\begin{proof}
In order to use theorem \ref{Prop(T)reduction} we will define a new simplicial complex $X'$ on which $\Gamma$ acts such that the conditions of theorem \ref{Prop(T)reduction} are fulfilled (in particular, the link of every vertex of $X'$ should be finite): \\
Let $r+1 \leq k \leq n-1$ such that $\lambda_k = \min_{\tau \in \Sigma (k-1)} \lambda (X_\tau) > C_k \dfrac{k}{k+1}$. Define the following weighted $2$-dimensional simplicial complex $X'$:
\begin{enumerate}
\item  $0$-simplices of $X'$ are $(k-1)$ (unordered) simplices of $X$. The weight of every such $0$-simplex $\sigma$ will be $m ' (\sigma) = m (\sigma)$ where $m$ is the weight given to $\sigma$ in $X$.
\item two different $0$-simplices in $X'$ are connected by an edge in $X'$ if they are both contained in a $k$-simplex in $X$, i.e. if $\lbrace u_0,...,u_{k-1} \rbrace, \lbrace v_0,...,v_{k-1} \rbrace$ are $(k-1)$-simplices in $X$, then they are connected by an edge in $X'$ if there is a $k$-simplex in $X$ which contain the vertices $u_0,...,u_{k-1}, v_0,...,v_{k-1}$. Note that a necessary condition for $\lbrace u_0,...,u_{k-1} \rbrace,  \lbrace v_0,...,v_{k-1} \rbrace$ to be connected by an edge is $\vert \lbrace u_0,...,u_{k-1} \rbrace \cap \lbrace v_0,...,v_{k-1} \rbrace \vert = k-1$. The weight of every such $1$-simplex $(\sigma_0, \sigma_1)$ will be $m ' ((\sigma_0,\sigma_1)) = m (\sigma_0 \cup \sigma_1)$ where $m$ is the weight given to $\sigma_0 \cup \sigma_1$ as a $k$-simplex in $X$.
\item There is a $2$-simplex in $X'$ in one of the following two cases: 
\begin{enumerate}
\item A $2$-simplex is defined by $3$ vertices of the form 
$$\lbrace u_0,...,u_{k-2},v \rbrace ,\lbrace u_0,...,u_{k-2},w \rbrace,\lbrace u_0,...,u_{k-2},x \rbrace $$ where $v,w,x$ are different from one another and $\lbrace u_0,...,u_{k-2},v,w,x \rbrace$ is a $(k+1)$-dimensional simplex in $X$. We denote the set of ordered simplices of this form as $\Sigma'^{(1)} (2)$. The weight of every such $2$-simplex $(\sigma_0, \sigma_1,\sigma_2) \in \Sigma'^{(1)} (2)$ will be $m ' ((\sigma_0,\sigma_1,\sigma_2)) = m (\sigma_0 \cup \sigma_1 \cup \sigma_2)$ where $m$ is the weight given to $\sigma_0 \cup \sigma_1 \cup \sigma_2$ as a $k+1$ simplex in $X$.
\item A $2$-simplex is defined by $3$ vertices of the form
$$\lbrace u_0,...,u_{k-3},v,w \rbrace, \lbrace u_0,...,u_{k-3},w,x \rbrace,\lbrace u_0,...,u_{k-3},v,x \rbrace$$ where $v,w,x$ are different from one another and $\lbrace u_0,...,u_{k-3},v,w,x \rbrace$ is a $k$-dimensional simplex in $X$. We denote the set of ordered simplices of this form as $\Sigma'^{(2)} (2)$. Note that for every $(\sigma_0, \sigma_1) \in \Sigma' (1)$ there are exactly $k-1$ different vertices $\sigma_2 \in \Sigma ' (0)$ such that $(\sigma_0, \sigma_1,\sigma_2) \in \Sigma'^{(2)} (2)$.  The weight of every such $2$-simplex $(\sigma_0, \sigma_1,\sigma_2) \in \Sigma'^{(2)} (2)$ will be $m ' ((\sigma_0,\sigma_1,\sigma_2))= \frac{\lambda_k m (\sigma_0 \cup \sigma_1 \cup \sigma_2)}{k}$ where $m$ is the weight given to $\sigma_0 \cup \sigma_1 \cup \sigma_2$ as a $k$ simplex in $X$.
\end{enumerate}    

\end{enumerate}
Since $X$ is connected with connected links, we get that $X'$ is connected with connected links and it is obvious that $\Gamma$ acts on $X'$ cocompactly and that the stabilizers of vertices are compact. For every fixed $(\sigma_0,\sigma_1) \in \Sigma' (1)$ we have the following:
$$ \sum_{(\sigma_0,\sigma_1,\sigma_2) \in \Sigma' (2)} m' ((\sigma_0,\sigma_1,\sigma_2)) =$$
$$ = \sum_{(\sigma_0,\sigma_1,\sigma_2) \in \Sigma'^{(1)} (2)} m (\sigma_0 \cup \sigma_1 \cup \sigma_2) +  \sum_{(\sigma_0,\sigma_1,\sigma_2) \in \Sigma'^{(2)} (2)} \dfrac{\lambda_k}{k} m (\sigma_0 \cup \sigma_1 \cup \sigma_2) =$$
$$ = C_k m (\sigma_0 \cup \sigma_1 ) + \dfrac{\lambda_k (k-1)}{k} m (\sigma_0 \cup \sigma_1 ) = ( C_k  + \dfrac{\lambda_k (k-1)}{k} ) m' ((\sigma_0,\sigma_1))$$
The reader should note that the above calculation doesn't include summing over different ordering - the vertices of $X'$ are unordered simplices and the summation is over triples of the form $(\sigma_0,\sigma_1,\sigma_2)$ where $\sigma_0,\sigma_1$ are fixed. This is the reason that there is no $(k+2)!$ multiplying $C_k$ in the last line. The above calculation yields $C'_1 = C_k  + \dfrac{\lambda_k (k-1)}{k}$ (again there is no $3!$ multiplying $C'_1$ because we don't sum over different ordering of each triple). By theorem \ref{Prop(T)reduction} we get that in order to show property (T), it will be enough to show that for every $\sigma \in \Sigma ' (0)$ we have 
$$\lambda (X'_\sigma) >  \dfrac{1}{2} (C_k  + \dfrac{\lambda_k (k-1)}{k})$$
We will do that by showing that $\lambda (X'_\sigma) \geq \lambda_k$. Indeed if this holds, then it is enough to show that:
$$\lambda_k >  \dfrac{1}{2} (C_k  + \dfrac{\lambda_k (k-1)}{k})$$
which is equivalent to 
$$\lambda_k >   C_k \dfrac{k}{k+1}$$ 
as assumed. \\ 
So we are left to show that for every $\sigma \in \Sigma ' (0)$ we have  we have that $\lambda (X'_\sigma) \geq \lambda_k$. Denote by $A_{X'_\sigma}$ the matrix corresponding to $\Delta^+$ on $X'_\sigma$ (as in lemma \ref{MatrixReduction} above) and by $A_{X_\sigma}$ the matrix corresponding to $\Delta^+$ on $X_\sigma$. So by lemma \ref{MatrixReduction}  we need to show that the smallest positive eigenvalue of $A_{X'_\sigma}$ is larger or equal to $\lambda_k$. Let us write $A_{X'_\sigma}$ explicitly:
$$A_{X'_\sigma} (\eta, \tau) = \begin{cases}   
								C_k  + \dfrac{\lambda_k (k-1)}{k} & \eta = \tau \\
								-\dfrac{m(\sigma \cup \eta \cup \tau)}{\sqrt{m(\sigma \cup \eta) m(\sigma \cup \tau)}} & (\sigma, \eta, \tau) \in \Sigma'^{(1)} (2) \\
								-\dfrac{\lambda_k}{k} & (\sigma, \eta, \tau) \in \Sigma'^{(2)} (2) \\
								0 & \text{otherwise}
								\end{cases}$$
(Note that in the cases $(\sigma, \eta, \tau) \in \Sigma'^{(2)} (2)$ we have that $m (\sigma \cup \eta \cup \tau)= m (\sigma \cup \eta)  = m (\sigma \cup \tau)$). \\
To complete the proof we need to recall some facts from matrix theory (see \cite{matrices} for further information and proofs): 
\begin{enumerate}
\item  The Kronecker product of matrices: for two matrices $M= (a_{ij})$ which is a $m \times m$ and $Q$ which is a $q \times q$ matrix the Kronecker product $M \otimes Q$ is a $mq \times mq$ matrix defined as
$$ M \otimes Q = \left( \begin{array}{cccc} a_{11} Q & a_{12} Q & ... & a_{1m} Q \\
											a_{21} Q & a_{22} Q & ... & a_{2m} Q \\    
 											\vdots & \vdots& \ddots & \vdots  \\
 											a_{m1} Q & a_{m2} Q & ... & a_{mm} Q 
 											\end{array} \right)$$
\item Denote by $I_l$ the $l \times l$ identity matrix. The Kronecker sum of $M$ and $Q$ (with the same dimensions as above) is defined to be
$$ M \oplus Q = I_q \otimes M + Q \otimes I_m$$
\item If $M$ is a $m \times m$ matrix with eigenvalues $\nu_1,...,\nu_m$ and $Q$ is a $q \times q$ matrix with eigenvalues $\mu_1,...,\mu_q$, then $M \oplus Q$ has the following $mq$ eigenvalues:
$$ \nu_1 + \mu_1,..., \nu_1 + \mu_q,\nu_2 + \mu_1,...,\nu_2 + \mu_q,..., \nu_m + \mu_q$$ 
\end{enumerate}
To finish the proof, it is enough to show that if we denote $B$ as the $k \times k$ matrix given by
$$ B =  \left( \begin{array}{cccc} \dfrac{\lambda_k (k-1)}{k} &-\dfrac{\lambda_k}{k} & ... & -\dfrac{\lambda_k}{k} \\
											-\dfrac{\lambda_k}{k} & \dfrac{\lambda_k (k-1)}{k} & ... & -\dfrac{\lambda_k}{k} \\    
 											\vdots & \vdots& \ddots & \vdots  \\
 											-\dfrac{\lambda_k}{k} & -\dfrac{\lambda_k}{k} & ... & \dfrac{\lambda_k (k-1)}{k}  
 											\end{array} \right)$$
Then
$$ A_{X'_\sigma} =  A_{X_\sigma} \oplus B$$
 
This is because $B$ has the eigenvalue $0$ with multiplicity $1$ and all its other eigenvalues are $\lambda_k$, and $A_{X_\sigma}$ has the eigenvalue $0$ with multiplicity $1$ and all its other eigenvalues are greater or equal to $\lambda_k$. Therefore $A_{X_\sigma} \oplus B$ has eigenvalue $0$ with multiplicity $1$ and all its other eigenvalues are greater or equal to $\lambda_k$ as needed. \\
To show that $A_{X'_\sigma} =  A_{X_\sigma} \oplus B$ holds recall that $A_{X_\sigma}$ can be written explicitly as a matrix with entries indexed by the vertices $u \in \Sigma_\sigma (0)$ as
$$A_{X_\sigma} (u, v) = \begin{cases}   
								C_{\sigma , 0}=C_k  & u = v \\
								\dfrac{-m_\sigma ((u,v))}{\sqrt{m_\sigma (u) m_\sigma (v)}} & (u,v) \in \Sigma_\sigma (1) \\
								0 & \text{otherwise}
					\end{cases} $$ 
So $I_k \otimes A_{X_\sigma}$ can be written as a matrix with entries indexed by couples of the form $(u,i)$ where $u \in \Sigma_\sigma (0)$ and $ 0 \leq i \leq k-1$ as
$$I_k \otimes A_{X_\sigma} ((u,i), (v,j)) = \begin{cases}   
								C_k  & u = v, i=j \\
								\dfrac{-m_\sigma ((u,v))}{\sqrt{m_\sigma (u) m_\sigma (v)}} & (u,v) \in \Sigma_\sigma (1), i=j \\
								0 & \text{otherwise}
					\end{cases} $$ 				
Fix an ordering of $\sigma$. A couple of the form $(u,i)$ can be identified as a $(k-1)$ unordered simplex in $X$ in the following way: identify $(u,i)$ with (the unordered) simplex $\sigma_i \cup u$ ($\sigma_i$ is the ordered simplex $\sigma$ with the $i$ entry omitted from it). Under this identification we get that $I_k \otimes A_{X_\sigma}$ is indexed by unordered $k-1$ simplices $\eta$ such that $(\sigma,\eta) \in \Sigma' (1)$ and
$$I_k \otimes A_{X_\sigma} (\eta, \tau) = \begin{cases}   
								C_k  & \eta=\tau \\
								-\dfrac{m(\sigma \cup \eta \cup \tau)}{\sqrt{m(\sigma \cup \eta) m(\sigma \cup \tau)}} & (\sigma, \eta, \tau) \in \Sigma'^{(1)} (2) \\
								0 & \text{otherwise}
					\end{cases} $$ 	 

The computation of $B \otimes I_{\vert \Sigma_\sigma (0) \vert}$ is similar - this matrix can be written as a matrix with entries are indexed by couples of the form $(u,i)$ where $u \in \Sigma_\sigma (0)$ and $ 0 \leq i \leq k-1$ as	
$$B \otimes I_{\vert \Sigma_\sigma (0) \vert} ((u,i),(v,j)) = \begin{cases}   
								 \dfrac{\lambda_k (k-1)}{k} & u=v, i=j \\
								-\dfrac{\lambda_k}{k} & u=v, i \neq j \\
								0 & \text{otherwise}
					\end{cases} $$ 
And under the same identification between couples $(u,i)$ and unordered $k-1$ simplices as before we get  
$$B \otimes I_{\vert \Sigma_\sigma (0) \vert} (\eta,\tau) = \begin{cases}   
								 \dfrac{\lambda_k (k-1)}{k} & \eta = \tau \\
								-\dfrac{\lambda_k}{k} & (\sigma, \eta, \tau) \in \Sigma'^{(2)} (2) \\
								0 & \text{otherwise}
					\end{cases} $$ 
So addition gives us the desired equality $$ A_{X'_\sigma} =  A_{X_\sigma} \oplus B$$
\end{proof}
We should remark that the above equality is nothing but a weighted version of a known fact about the Laplacian matrix of the Cartesian product of two graphs. Given two graphs $(V_1,E_1),(V_2,E_2)$, the Cartesian product $(V_1,E_1) \times (V_2,E_2)$ is the graph given by $V = V_1 \times V_2$ and $E$ defined by $(v_1,v_2) \sim (v_1',v_2')$ if $v_1=v_1'$ and $v_2 \sim v_2'$ or $v_1 \sim v_1'$ and $v_2 = v_2'$. It is known (see \cite{Fiedler}[Theorem 3.4]) that the Laplacian matrix of the Cartesian product $A_{(V,E)}$ is equal to $A_{(V_1,E_1)} \oplus A_{(V_2,E_2)}$. In our case, as a graph $X_\sigma '$ is the Cartesian product of $X_\sigma$ with the complete graph of $k$ vertices. However, the use of weights in the Laplacian prevented us to use \cite{Fiedler}[Theorem 3.4] in its original form and required us to prove the above version of it in our weighted situation.
\bibliographystyle{alpha}
\bibliography{bibl}

\end{document}